\documentclass[11pt]{amsart}
\usepackage{amssymb,amsmath}
\usepackage{amsthm}
\usepackage{amsfonts}
\usepackage[dvips]{graphicx}
\usepackage{enumerate}

\newcommand{\Li}{\mathcal{L}}
\newcommand{\T}{\mathbb{T}}
\newcommand{\C}{\mathcal{C}}
\newcommand{\R}{\mathbb{R}}
\newcommand{\Z}{\mathbb{Z}}

\renewcommand{\epsilon}{\varepsilon}

\def\ep{\noindent{\hfill $\fbox{\,}$}\bigskip\newline}
\newtheorem{problem}{\sc Problem}
\newtheorem{conjecture}{\sc Conjecture}
\newtheorem{question}{\sc Question}

\newtheorem{theorem}{Theorem}[section]

\newtheorem{cor}[theorem]{Corollary}

\newtheorem{prop}[theorem]{Proposition}

\title{Cohomology free systems and the first Betti number.}
\author{Federico Rodriguez Hertz}
\author{Jana Rodriguez Hertz}
\address{}
\email{jana@fing.edu.uy}
\email{frhertz@fing.edu.uy}

\date{\today}
\begin{document}

\begin{abstract}
We prove that a cohomology free flow on a manifold $M$ fibers over a diophantine translation on $\T^{\beta_1}$ where $\beta_1$ is the first Betti number of $M$.
\end{abstract}

\maketitle
\begin{section}{Introduction}
Let $X$ be a vector field over a manifold M. A
smooth\footnote{smooth means $C^{\infty}$ along this paper}
function $\phi$ (also called a cocycle) is called a coboundary if
there is a smooth function $u$ (called a transfer function)
solving the equation:
$$
\Li_Xu=\phi
$$
where $\Li_Xu$ is the Lie derivative of $u$ along $X$. Two
cocycles are called cohomologous if their difference is a
coboundary. Notice that for a cocycle to be a coboundary it is
necessary to have zero average with respect to any measure
invariant under the flow generated by $X$. We say that $X$ is a
cohomology free vector field if every cocycle is cohomologous to a
constant cocycle. There are other definitions considering finite
differentiability (see for instance \cite{kr}).

Let us give the main examples of a cohomolgy free vector field. A
vector $\alpha\in\R^N$ is said to be diophantine if there are
constants $c>0$, $\gamma>0$ such that
$$
|n\cdot\alpha|\geq\frac{c}{|n|^{\gamma}}
$$
for every $n\in\Z^N\setminus\{0\}$. A vector $\alpha$ in $\R^N$
defines a constant vector field on $\T^N=\R^N/\Z^N$. It is known
that $X_{\alpha}\equiv\alpha$ in $\T^N$ is cohomology free if and
only if $\alpha$ is diophantine. Note that these are the only
known examples by now. And in fact, the following conjecture was
posed \cite{h}, \cite{kr}:
\begin{conjecture}
Any cohomology free vector field is smoothly conjugated to a diophantine vector field on a torus.
\end{conjecture}
The main theorem is an approach to this conjecture.  Along this
paper $M$ will be a compact and connected manifold.
\begin{theorem}\label{teo}
Let $X$ be a cohomology free vector field on $M$. Then there is a
fibration $p:M\to\T^{\beta_1}$, where $\beta_1$ is the first Betti
number of $M$, such that $p_*X=\alpha$ is a constant diophantine
vector field.
\end{theorem}
As a corollary we get the following:
\begin{cor}
Let $X$ be a cohomology free vector field on $M$. Then the first Betti number of $M$ is less than or equal to the
dimension of $M$. Moreover, if the first Betti number of $M$ is equal to the dimension of $M$ then $M$ is a torus
and the vector field is smoothly conjugated to a diophantine vector field.
\end{cor}
To prove the theorem, we shall prove the following proposition which we find interesting in itself.
\begin{prop}\label{invfib}
Let $X$ be a cohomology free vector field on $M$, and assume there
is a fibration $p:M\to N$ and a vector field $Y$ in $N$ such that
$p_*X=Y$, then $Y$ is also cohomology free.
\end{prop}
Let us point out that in \cite{ls} the authors get similar
results. The main novelty here is that are only assumping the
underlying manifold is compact and connected.\par
In the final section we shall discuss some related problems and
results.\newline
\noindent {\bf Acknowledgment:} We want to thank Giovanni Forni
for the reference \cite{gw}. Also, we are grateful to the referees
for pointing out that the proof of theorem \ref{teo} needed some
extra explanation. The first author thanks the organizers of the
conference for all the support and hospitality.
\end{section}

\begin{section}{Proof of proposition \ref{invfib}}
Let us begin with some known fact:
\begin{prop}\label{unique}
If $X$ is a cohomology free vector field, then it is uniquely ergodic and its unique invariant probability measure
is a smooth volume form. Therefore it is minimal.
\end{prop}

The reader can find a proof of this proposition in \cite{kr} page
21. Let us notice that it is also true that this volume form is
the only invariant distribution.\par

{\em Proof of proposition \ref{invfib}.} Take $p:M\to N$ with
$p_*X=Y$ as in the statement, and let us prove that $Y$ is
cohomology free. Take a function $\psi:N\to\R$ and define
$\phi:M\to\R$ by $\phi=\psi\circ p$. By hypothesis, there is a
function $\hat{u}:M\to\R$ and $c\in\R$ such that
$\Li_X\hat{u}=\phi+c$. By proposition \ref{unique}, we have a
volume form $\omega$ invariant by $X$. Let us denote $\omega$
restricted to $M_y=p^{-1}(y)$ by $\omega_y$, $y\in N$, and define
$u:N\to\R$ by
$$
u(y)=\int_{M_y}\hat{u}\omega_y
$$
As $\omega$ is an invariant form and $p_*X=Y$, a straight forward computation gives
$$
\Li_Yu(y)=\int_{M_y}\bigl(\Li_X\hat{u}\bigr)\omega_y=\int_{M_y}\bigl(\phi+c\bigr)\omega_y
$$
As $\phi$ is constant on each fiber, we get that
$$
\Li_Yu(y)=\bigr(\psi(y)+c\bigl)\mu_y(M_y)
$$
where $\mu_y$ is the measure generated by $\omega_y$. Now,
invariance of $\omega$ implies invariance of
$x\to\mu_{p(x)}(M_{p(x)})$ and, by minimality, it follows
$\mu_y(M_y)=1$ for any $y\in N$, what finishes the proof.\ep
\end{section}

\begin{section}{Proof of theorem \ref{teo}}
Take a smooth map $\hat{p}:M\to\T^{\beta_1}$ realizing the
Hurewicz homomorphism
$$
\hat{p}^{\#}:\pi_1(M)\to H^1(M)=\Z^{\beta_1}
$$
We have that $d_x\hat{p}(X(x))$ is a smooth map from $M$ to $\R^{\beta_1}$ and hence there is a map $u:M\to\R^{\beta_1}$
and a vector $\alpha\in\R^{\beta_1}$ such that
$$
\Li_Xu(x)=-d_x\hat{p}(X(x))+\alpha
$$
where $\Li$ acts componentwise. Notice that we can define $\hat{u}:M\to\T^{\beta_1}$ just projecting $u$ into the torus.
Define now $p:M\to\T^{\beta_1}$ by $p=\hat{p}+\hat{u}$. Then it follows that
$$
p_*X(x)=d_x p(X(x))=d_x\hat{p}(X(x))+d_x\hat{u}(X(x))=d_x\hat{p}(X(x))+\Li_Xu(x)=\alpha
$$
This proves that $p_*X=\alpha$. In order to see that $p$ is a
fibration, we only have to prove that $p$ is a
submersion. As $u+const$ is also a transfer function, we may
assume without loss of generality that $p(z_0)=0$ for some point
$z_0\in M$. As $X$ generates a minimal flow the image of $p$ will
equal the closure of the orbit of $0$ under the linear flow
generated by $\alpha$. So the image of $p$ is a compact connected
subgroup of $\T^{\beta_1}$ and thus it is a rational subtorus. On
the other hand, as $p=\hat{p}+\hat{u}$, $p$ is homotopic to
$\hat{p}$ and thus $p^{\#}=\hat{p}^{\#}$ is the Hurewicz
homomorphism. Finally, as the Hurewicz homomorphism is onto,
putting all together we get that the image of $p$ is a rational
subtorus of $\T^{\beta_1}$ with $\pi_1$ equals to $\Z^{\beta_1}$
wich must be the whole $\T^{\beta_1}$. So we get that $p$ is onto.
Now, let us call
$$
R=\{x\in M\;\;\mbox{such that the rank of}\;\;
d_xp\;\;\mbox{is}\;\;\beta_1\}
$$
$R$ is open because having maximal rank is an open condition. As
$p_*X$ is contant, it is straightforward that $R$ is invariant
under the flow generated by $X$. Again, as this flow is minimal,
either $R$ is empty or is the whole $M$. Since $p$ is surjective,
it follows, by using Sard's theorem, that $R$ is not empty; and
thus $R$ must be $M$ and $p$ is a submersion. This shows $p$ is a
fibration. \par%
 Now proposition \ref{invfib} implies that the constant vector field
$\alpha$ on $\T^{\beta_1}$ is cohomology free. Hence the results
in \cite{kr} page 19 imply that $\alpha$ is diophantine, and
theorem \ref{teo} follows.
\end{section}

\begin{section}{Final remarks and open questions}
It is worth noting that an equivalent statement of our result in
the setting of diffeomorphisms would be:
\begin{theorem}\label{difeos}
Let $f:M\to M$ be a cohomology free diffeomorphism. Then there is a fibration $p:M\to \T^{\hat b_1}$ where $\hat b_1=\dim \ker (f_*-Id)$ and $f_*$ is the induced action on  $H^1(M,\R)$. Moreover $p\circ f= p+\alpha$ with $\alpha$ diophantine.
\end{theorem}
The proof of theorem above follows by either restating the proof
of theorem \ref{teo} for diffeomorphisms, or by taking the
suspension flow of $f$ and getting a cohomology free vector field
on a manifold $\hat M$. The key property to be used is the
following formula:
$$
\hat b_1+1=\mbox{first Betti number of } \hat M
$$
In \cite{gw} the reader may find two related problems. A globally
hypoelliptic vector field $X$ is one verifying that all solutions
$u$ (in the distribution sense) of the following equation:
$$
\Li_Xu=\phi
$$
with $\phi$ in $\C^{\infty}$ are, in fact, $\C^{\infty}$
functions. Clearly a cohomology free vector field is globally
hypoelliptic. It seems likely that globally hypoelliptic and
cohomology free vector fields are the same thing (see \cite{cc}).
In \cite{gw} the authors posed the following questions:
\begin{question}
Is any globally hypoelliptic vector field smoothly conjugated to a diophantine vector field on a torus?
\end{question}
\begin{question}
What about a globally hypoelliptic homogeneous flow?
\end{question}
For these problems we propose the following approaches:
\begin{problem}
An homogeneous flow on a nilmanifold other than a torus has an invariant distribution other than multiples of volume. (It seems that Flaminio and Forni have some results related to this).
\end{problem}
\begin{problem}
Given a globally hypoelliptic vector field $X$, there is a fibration $p:M\to\T^{\beta_1}$ where $\beta_1$ is the first Betti number of $M$ such that $p_*X=\alpha$ is a constant diophantine vector field and the fibers of the fibration are simply connected.
\end{problem}
We want to finish with the following:
\begin{problem}
Find a uniquely ergodic $\C^{\infty}$ system with positive entropy.
\end{problem}
\end{section}

\bibliographystyle{alpha}

\end{document}